\documentclass[12pt]{article}
\usepackage[T2A]{fontenc}
\usepackage[cp1251]{inputenc}
\usepackage[russian]{babel}
\usepackage{amssymb,amsmath,latexsym,amsthm}
\usepackage{amsbsy,amsfonts,amssymb,mathrsfs,amscd}
\usepackage{eucal}
\usepackage{graphicx}
\usepackage{indentfirst}
\usepackage{verbatim}
\usepackage{enumerate}
\usepackage[pdftex,unicode,colorlinks,linkcolor=blue,
citecolor=red,bookmarksopen,pdfhighlight=/N]{hyperref}

\theoremstyle{plain}

\newtheorem{lemma}{Лемма}
\newtheorem{theorem}{Теорема}
\newtheorem{corollary}{Следствие}
\newtheorem{proposition}{Предложение}

\newtheorem*{thA}{Теорема A}
\newtheorem*{thB}{Теорема B}

\newtheorem*{thD}{Теорема C}
\newtheorem*{thHL}{Теорема Адамара\,--\,Линделёфа}

\theoremstyle{definition}

\newtheorem{remark}{Замечание}

\renewcommand{\leq}{\leqslant}
\renewcommand{\geq}{\geqslant}

\newcommand{\rad}{\text{\tiny\rm rad}}
\newcommand{\RR}{\mathbb{R}} 
\newcommand{\CC}{\mathbb{C}} 
\newcommand{\NN}{\mathbb{N}}

\DeclareMathOperator{\mes}{mes}
\DeclareMathOperator{\dens}{dns} 
\DeclareMathOperator{\Zero}{\text{\sf Zero}}

\DeclareMathOperator{\type}{\text{\sf type}} 
\DeclareMathOperator{\Hol}{Hol}

\DeclareMathOperator{\rh}{\text{\rm rh}}
 \DeclareMathOperator{\lh}{\text{\rm lh}}

\DeclareMathOperator{\dd}{\,{\mathrm d\!}}
\renewcommand{\Re}{{\rm Re \,}}
\renewcommand{\Im}{{\rm Im \,}}

\numberwithin{equation}{section}

\title{Рост целых функций экспоненциального типа и характеристики распределений точек вдоль прямой}
\author{А.\,Е.~Салимова, Б.\,Н.~Хабибуллина}
\date{4 мая 2021 г.}
\begin{document}
\maketitle

\maketitle {
\begin{quote}
\noindent{\bf Аннотация.} Для пары распределений точек ${\sf Z}$ и ${\sf W}$ конечной верхней плотности на комплексной плоскости $\mathbb C$ с вещественной осью $\mathbb R$ даются несколько версий необходимых и одновременно  достаточных условий на их расположение, при которых для любой целой функции экспоненциального типа $g\neq 0$, обращающейся в нуль на ${\sf W}$, либо существует целая функция экспоненциального типа $f\neq 0$, обращающаяся в нуль на ${\sf Z}$ и удовлетворяющая одному из двух вариантов ограничений: 
\begin{itemize}
\item $\bigl|f(iy)\bigr|\leq \bigl|g(iy)\bigr|$ для всех $y\in \mathbb R$, т.е. всюду на мнимой оси $i\mathbb R$,
\item $\ln \bigl|f(iy)\bigr|\leq \ln \bigl|g(iy)\bigr|+o\bigl(|y|\bigr)$ при $y\to \pm \infty$,
\end{itemize}
либо для любого числа $\varepsilon >0$ найдётся целая функция экспоненциального типа $f\neq 0$, обращающаяся в нуль на ${\sf Z}$ и удовлетворяющая неравенству $\ln \bigl|f(iy)\bigr|\leq \ln \bigl|g(iy)\bigr|+\varepsilon |y|$ при всех $y\in \mathbb R\setminus E$, где $E\subset \mathbb R$ --- множество конечной линейной лебеговой  меры.
Исследование проведено в рамках обобщения и развития классической теоремы П.~Мальявена и Л.~А.~Рубела 1960-х гг., в которой  был рассмотрен только случай расположения  ${\sf Z}\subset  \RR^+$ и ${\sf W}\subset  \RR^+$	на положительной полуоси $\mathbb R^+\subset \mathbb R$.  Критерии даются в терминах мажорирования логарифмических характеристик и (суб)мер для   ${\sf Z}$, соответствующими  
логарифмическими характеристиками и (суб)мерами для  ${\sf W}$. При этом в последнем третьем варианте никаких дополнительных требований на ${\sf Z}$ и ${\sf W}$ не накладывается, а в первом и втором вариантах предполагается  асимптотическая отделённость углами от мнимой оси для ${\sf Z}$ и ${\sf W}$ и либо расположение ${\sf W}$ полностью в правой или левой  полуплоскости, либо  условие типа Линделёфа для ${\sf W}$ вдоль мнимой оси $i\mathbb R$ об определённой симметричности мнимых частей обратных величин $1/{\sf w}$ к ${\sf w}\in {\sf W}$ вида 
$\Bigl|\sum_{1\leq |{\sf w}|\leq r}\Im (1/{\sf w}) \Bigr|=O(1)$ при $r\to +\infty$.
\medskip

\noindent{\bf Ключевые слова:} {целая функция экспоненциального типа, распределение корней, рост целой функции, логарифмические характеристики и меры, условие Линделёфа}


\end{quote}
}

\section{Введение}\label{s10}
\subsection{Обозначения и соглашения}\label{11def} 
Одноточечные множества $\{x\}$ часто записываем без фигурных скобок, т.е. просто как $x$. 
Как обычно, $\mathbb N:=\{1,2, \dots\}$ --- множество 
{\it натуральных чисел,\/}  $\NN_0:=0\cup \NN$. 
Множество $\RR$ {\it действительных чисел,\/} со стандартными порядковой
($\le$, $\sup/\inf$), алгебраической и топологической структурами в основном рассматривается  как  {\it вещественная ось\/} в  {\it комплексной плоскости\/} $\CC$; $i\RR$ --- {\it мнимая ось},  
$\overline \RR:=-\infty\cup \RR\cup +\infty$ --- {\it расширенная действительная прямая\/} с двумя концами $\pm\infty\notin \RR$, дополненная неравенствами $-\infty\leq x\leq +\infty$ для любого $x\in \overline \RR$ и снабжённая естественной порядковой топологией, а $\overline \NN:=\NN\cup +\infty$, $\overline \NN_0:=\NN_0\cup +\infty$.
{\it Интервалы на\/} $\overline{\RR}$ --- связные подмножества в $\overline \RR$, такие, как  {\it отрезок\/} 
$[a,b]:=\{x\in \overline\RR\colon a\leq x\leq b\}$ с {\it концами\/}  $a,b\in \overline \RR$, где   $[a,b]=\varnothing$ --- {\it пустое множество\/} при  $a< b$, а также $(a,b]:=[a,b]\setminus a$, $[a,b):=[a,b]\setminus b$ и {\it открытый интервал\/}
$(a,b):=(a,b]\cap [a,b)$.  По определению $\inf \varnothing :=+\infty$ и  $\sup \varnothing :=-\infty$.
{\it Правые и левые открытые полуплоскости\/} обозначаем  соответственно как
$\CC_{\rh}:= \{z\in \CC \colon \Re z>0\}$ и  ${\CC_{\lh}}:=-\CC_{\rh}$. 
Через  $\overline D(r):=\bigl\{z\in \CC\colon |z|\leq r\bigr\}$ обозначаем {\it замкнутый круг радиуса $r\in \RR^+$ с центром в нуле}. Для $x\in X\subset \overline \RR$ полагаем 
$x^+:=\sup\{0,x \}$, а 
{\it расширенной числовой функции\/} $f\colon S\to \overline \RR$ сопоставляем её {\it положительную часть\/} $f^+\colon s\underset{\text{\tiny $s\!\in\! S$}}{\longmapsto} (f(s))^+\in \overline{\RR}^+$ 
 и {\it отрицательную часть\/} $f^-:=(-f)^+\colon S\to \overline{\RR}^+ $.
Если в записи суммы 
верхний предел суммирования меньше нижнего предела  суммирования или суммирование ведётся по пустому множеству, то по определению  сумма  равна нулю.

\subsection{Распределения точек на комплексной плоскости}\label{subs1_1}  
 Каждому {\it распределению точек\/}  ${\sf Z}=\{{\sf z}_{\tt j}\}$ из $\CC$, состоящему из пронумерованных не более чем счётным количеством   индексов ${\tt j}$ точек ${\sf z}_{\tt j}\in \CC$, сопоставляем {\it считающую меру\/}  
 \cite[0.1.2]{Khsur}
\begin{subequations}\label{n}
\begin{align}
n_{\sf Z}&\colon  S\underset{S\subset \CC}{\longmapsto} 
\sum_{{\sf z}_{\tt j}\in S}1\in   \overline \NN_0
\tag{\ref{n}n}\label{df:divmn}\\
\intertext{--- число точек ${\sf z}_{\tt j}$, попавших в $S$. При этом то же обозначение}
n_{\sf Z}&\colon z\underset{\text{\tiny $z\in \CC$}}{\longmapsto} n_{\sf Z}(z)
=\sum_{{\sf z}_{\tt j}=z}1\in  \overline \NN_0
\tag{\ref{n}z}\label{df:divz}\\
\intertext{используется и для {\it считающей функции\/} распределения точек  ${\sf Z}$, а}
n_{\sf Z}^{\rad}(r)&\overset{\eqref{df:divz}}{\underset{r\in \RR^+}{:=}}
\sum_{|z|\leq r}n_{\sf Z}(z)=
\sum_{|{\sf z}_{\tt j}|\leq r}1
\overset{\eqref{df:divmn}}{=}n_{\sf Z}\bigl( \overline D(r)\bigr)\in   \overline \NN_0
\tag{\ref{n}r}\label{nrad}
\end{align}
\end{subequations}
--- 
{\it радиальная считающая функция\/} для ${\sf Z}$ по {\it замкнутым кругам\/} 
$\overline D(r)$.

Распределения точек    ${\sf Z}$ и  ${\sf Z'}$  {\it совпадают,\/} или {\it равны\/,} и пишем  ${\sf Z}={\sf Z'}$  если у них одни и те же считающие меры или функции $n_{{\sf Z}}\overset{\eqref{df:divmn}}{=}n_{\sf Z'}$, а  включение ${\sf Z}\subset {\sf Z}'$ означает, что  выполнено неравенство $n_{\sf Z}\leq n_{\sf Z'}$. {\it Объединение\/}  ${\sf Z}\cup {\sf Z'}$ определяется считающей мерой или функцией  $n_{{\sf Z}\cup {\sf Z'}}:=n_{\sf Z}+n_{\sf Z'}$, а при    ${\sf Z}\subset  {\sf Z'}$ {\it разность\/} ${\sf Z'}\setminus  {\sf Z}$ определяется считающей мерой или функцией  $n_{{\sf Z'}\setminus {\sf Z}}:=n_{\sf Z'}-n_{\sf Z}$.  Точка $z\in \CC$ {\it принадлежит\/} ${\sf Z}$, т.е. $z\in {\sf Z}$, если $n_{\sf Z}(z)> 0$ для считающей функции \eqref{df:divz}.  

{\it Верхняя плотность\/}  распределения точек ${\sf Z}$ (при порядке $1$) определяется как \cite{Boas}, \cite{Levin56}, \cite{Levin96}
\begin{equation}\label{dZ}
\overline \dens_{\sf Z}:=\limsup_{r\to +\infty} \frac{n_{\sf Z}^{\rad}(r)}{r}\in \overline \RR^+. 
\end{equation}

Ключевую роль далее будут играть логарифмические характеристики для распределений точек на $\CC$, введённые для {\it положительных\/} распределений точек ${\sf Z}\subset \RR^+$ в основополагающей для настоящего исследования 
статье   П. Мальявена и Л. А. Рубела \cite{MR} (см. также монографию  Л. А. Рубела \cite{RC}) и распространённые на произвольные {\it комплексные\/} распределения точек   ${\sf Z}\subset \CC$ в работах Б. Н. Хабибуллина 
\cite{KhaD88}, \cite{Kha89}, \cite{kh91AA}, \cite[3.2]{Khsur}, \cite{Kha01l}. 

Определим {\it  правый и левый характеристические логарифмы} для ${\sf Z}\subset \CC$ как
\begin{subequations}\label{logZC}
\begin{align}
l_{{\sf Z}}^{\rh }(r)&:=\sum_{\substack{0 < |{\sf z}_k|\leq r\\{\sf z}_k \in \CC_{\rh }}} \Re \frac{1}{{\sf z}_k}
=\sum_{0< |{\sf z}_k|\leq r} \Re^+ \frac{1}{{\sf z}_k}, \quad 0<r\leq +\infty, 
\tag{\ref{logZC}r}\label{df:dD+}\\
l_{{\sf Z}}^{\lh }(r)&:=
\sum_{\substack{0< |{\sf z}_k|\leq r\\{\sf z}_k \in \CC_{\lh }}} 
-\Re \frac{1}{{\sf z}_k}=\sum_{0< |{\sf z}_k|\leq r} \Re^- \frac{1}{{\sf z}_k}, \quad 0<r\leq +\infty,
\tag{\ref{logZC}l}\label{df:dD-}
\end{align}
\end{subequations}
а также  {\it правую\/} и {\it левую  логарифмические меры интервалов\/}    
\begin{subequations}\label{df:l}
\begin{align}
l_{{\sf Z}}^{\rh }(r, R)&\overset{\eqref{df:dD+}}{:=}
l_{{\sf Z}}^{\rh }(R)-l_{{\sf Z}}^{\rh }(r),
\quad 0< r < R \leq +\infty,
\tag{\ref{df:l}r}\label{df:dDl+}\\
l_{{\sf Z}}^{\lh }(r, R)&\overset{\eqref{df:dD-}}{:=}l_{{\sf Z}}^{\lh }(R)-l_{{\sf Z}}^{\lh }(r), \quad 0< r < R \leq +\infty,
\tag{\ref{df:l}l}\label{df:dDl-}
\\
\intertext{которые порождают {\it логарифмическую субмеру  интервалов\/} для  $\sf Z$:} 
l_{{\sf Z}}(r, R)&:=\max \{ l_{{\sf Z}}^{\lh }(r, R), l_{{\sf Z}}^{\rh }(r,
R)\}, \quad 0< r < R \leq +\infty ,
\tag{\ref{df:l}m}\label{df:dDlL}
\end{align}
\end{subequations}
где для ${\sf Z}=\varnothing$ по определению $l_{\varnothing}(r,R)\equiv 0$ при всех $0< r < R \leq +\infty$.

\subsection{Целые функции, предшествующие результаты и постановки задач} 
Кольцо $\Hol (\CC)$ над $\CC$ состоит из всех голоморфных функций на $\CC$, т.е. $\Hol (\CC)$ --- кольцо {\it целых функций.\/} Через $\Hol_*(\CC):=\bigl\{f\in \Hol(\CC)\colon f\neq 0\bigr\}$ обозначаем множество  всех ненулевых целых функций. Через $\Zero_f$ обозначаем {\it распределение всех  корней\/} целой функции $f\neq 0$ со считающей функцией $n_{\Zero_f}$ в смысле  \eqref{df:divz}, равной  в каждой точке $z\in \CC$ кратности корня функции $f$ в точке $z$.
Целая функция $f\neq 0$ на  $\CC$ \textit{обращается в нуль\/} на распределении  точек ${\sf Z}$ и пишем $f({\sf Z})=0$, если ${\sf Z}\subset {\Zero}_f$.

Целую функцию $f$ называют \textit{целой функции экспоненциального типа} (пишем {\it ц.ф.э.т}), 
если конечен её тип (при порядке $1$) \cite{Boas}, \cite{Levin56}, \cite{Levin96}, \cite[2.1]{KhI}
 \begin{equation}\label{typef}
\type_f:=\limsup_{z\to \infty}\frac{\ln |f(z)|}{|z|}.
\end{equation} 

Распределение точек   ${\sf Z}=\{{\sf z}_{\tt j}\}\subset \CC$ 
порождает идеал \cite{MR},  \cite[гл. 22]{RC}, \cite{SalKha20}
\begin{equation*}
I({\sf Z}):=\bigl\{f\in \Hol (\CC)\colon f({\sf Z})=0\bigr\}\subset \Hol (\CC)
\end{equation*} 
в кольце $\Hol (\CC)$, а также идеал в кольце всех ц.ф.э.т.
\footnote{В \cite{MR} и  \cite{RC} идеал $I^1({\sf Z})$ обозначен соответственно  как $\mathcal F ({\sf Z})$ и $F({\sf Z})$.}
\begin{equation*}
I^1({\sf Z}):=I({\sf Z})\cap \bigl\{f\in \Hol (\CC)\colon {\type}_f\overset{\eqref{typef}}{<}+\infty \bigr\},
\end{equation*} 
для которых  полагаем 
\begin{equation*}
I_*({\sf Z}):=I({\sf Z})\cap \Hol_*(\CC),
\quad 
I_*^1({\sf Z}):=I^1({\sf Z})\cap \Hol_*(\CC).
\end{equation*}

Всюду далее ${\sf Z}\subset \CC$,  ${\sf W}=\{{\sf w}_{\tt j}\}\subset \CC$ --- {\it распределения точек конечной верхней плотности\/} 
 \begin{equation}\label{dnsZW}
\overline\dens_{\sf Z}+\overline\dens_{\sf W}<+\infty.
\end{equation}

Через  $\mes$ обозначаем линейную лебегову меру евклидовой длины на $\RR$.

\begin{thA}[{\cite[основная теорема]{KhaD88}, \cite[основная теорема]{Kha89}, \cite[теорема 3.2.1]{Khsur}}] Пусть при \eqref{dnsZW}  ещё и   ${\sf W}\subset \CC_{\rh}$ из правой полуплоскости. Тогда эквивалентны три утверждения: 
\begin{enumerate}[{\rm I.}]
\item\label{fgie} Для  любой функции $g\in I_*^1({\sf W})$ при любом $\varepsilon \in \RR^+\setminus 0$
найдутся функция $f\in I^1_*({\sf Z})$ и  борелевское множество $E\subset \RR$, для которых 
\begin{equation}\label{fgiRe}
\ln \bigl|f(iy)\bigr|\leq \ln \bigl|g(iy)\bigr|+\varepsilon |y|\quad \text{при всех  $y\in \RR\setminus E$  и   $\mes E<+\infty$}.
\end{equation}
\item\label{fgiiie} При любом $\varepsilon \in \RR^+\setminus 0$
существуют пара $f\in I_*^1({\sf Z})$ и  $g\in I^1_*({\sf W})$ c $\Zero_g \cap \CC_{\rh} ={\sf W}$ и  борелевское множество $E\subset \RR$, удовлетворяющие   \eqref{fgiRe}. 
\item\label{fgiie} Для любого $\varepsilon \in \RR^+\setminus 0$ найдётся $C\in \RR^+$, для которого 
\begin{equation}\label{Zlde}
l_{\sf Z}(r,R)\leq l_{\sf W}(r,R) +\varepsilon \ln \frac{R}{r}
+C\quad \text{при всех\/ $0< r<R<+\infty$}.
\end{equation}
\end{enumerate}
\end{thA}

Распределение точек  ${\sf Z} =\{ {\sf z}_{\tt j}\} \subset \CC$  {\it асимптотически отделено углами  от\/}
$i\RR$,  если \cite[(1.2)]{SalKha20}
\begin{equation}\label{con:dis}
\left(\liminf_{{\tt j}\to\infty}\frac{\bigl| \Re {\sf z}_{\tt j} \bigr|}{ |{\sf z}_{\tt j} |} >0 \right)
\underset{\text{\it или}}{\Longleftrightarrow}
\left(\limsup_{{\tt j}\to\infty}\frac{\bigl| \Im {\sf z}_{\tt j} \bigr|}{ |{\sf z}_{\tt j} |} <1 \right).
\end{equation}
 Пара эквивалентных ограничений \eqref{con:dis} геометрически означает, что  найдётся  пара непустых открытых вертикальных углов, содержащих $i\RR\setminus 0$, для которой точки ${\sf z}_{\tt j}$ лежат вне этой пары углов  при всех ${\tt j}$ за исключением их конечного числа.

\begin{thB}[{\cite[основная теорема]{SalKha20}}] Пусть при \eqref{dnsZW}  и ${\sf Z}\subset \CC$, и  ${\sf W}\subset \CC_{\rh}$  асимптотически отделены углами от  $i\RR$ в смысле \eqref{con:dis}. 
Тогда эквивалентны три утверждения: 
\begin{enumerate}[{\rm I.}]
\item\label{fgi} Для любой функции $g\in I_*^1({\sf W})$ с распределением корней, асимптотически отделённым углами  от $i\RR$,  найдётся функция $f\in I^1_*({\sf Z})$ с ограничением 
\begin{equation}\label{fgiR}
\ln \bigl|f(iy)\bigr|\leq \ln \bigl|g(iy)\bigr|+o(|y|)\quad \text{при  $y\to  +\infty$}.
\end{equation}
\item\label{fgiii} Существует пара $f\in I_*^1({\sf Z})$ и  $g\in I^1_*({\sf W})$ c $\Zero_g \cap \CC_{\rh} ={\sf W}$ удовлетворяющая  \eqref{fgiR}. 
\item\label{fgii} Существуют $C\in \RR^+$ и ограниченная  функция 
$d\colon \RR^+\to \RR^+$ с $d(x)\underset{x\to +\infty}{=}o(1)$, для которых выполнены неравенства 
\begin{equation}\label{Zld}
l_{\sf Z}(r,R)\leq l_{\sf W}(r,R) +d(R)\ln \frac{R}{r}
+C\quad \text{при всех\/ $0< r<R<+\infty$}.
\end{equation}
\end{enumerate}
\end{thB}

\begin{thD}[{\cite[теорема 1]{SalKha21}}] Пусть 
при \eqref{dnsZW} распределения точек ${\sf Z}\subset \CC$ и ${\sf W}\subset \CC_{\rh}$, 
 асимптотически отделённы углами от  $i\RR$  в смысле \eqref{con:dis}. 
Тогда эквивалентны три утверждения: 
\begin{enumerate}[{\rm I.}]
\item\label{fgi0} Для любой функции $g\in I_*^1({\sf W})$  найдётся функция $f\in I^1_*({\sf Z})$ с ограничением 
\begin{equation}\label{fgiR0}
 \bigl|f(iy)\bigr|\leq  \bigl|g(iy)\bigr|\quad \text{при всех  $y\in \RR$}.
\end{equation}
\item\label{fgiii0} Существует пара $f\in I_*^1({\sf Z})$ и  $g\in I^1_*({\sf W})$ c $\Zero_g \cap \CC_{\rh} ={\sf W}$, удовлетворяющая   \eqref{fgiR0}. 
\item\label{fgii0} Существуют $C\in \RR^+$, для которого  
\begin{equation}\label{Zld0}
l_{\sf Z}(r,R)\leq l_{\sf W}(r,R) +C\quad \text{при всех\/ $0< r<R<+\infty$}.
\end{equation}
\end{enumerate}
\end{thD}

В нашей статье выполняются следующие две задачи.
\begin{enumerate}
\item Полное удаление условия  о  расположении  ${\sf W}$ в правой полуплоскости из теоремы A  при  некотором изменении  утверждения \ref{fgiiie}  в  теореме \ref{thAW} из \S~\ref{S3}.
\item Замена  расположения  ${\sf W}$ в правой полуплоскости $\CC_{\rh}$ из  теорем B и  C  на условие 
\begin{equation}\label{Lind}
\sup_{\tt j}\biggl|\sum_{1<|{\sf w}_{\tt j}|\leq r} \Im \frac{1}{{\sf w}_{\tt j}}\biggr|<+\infty
\end{equation}
в теоремах \ref{thBWB} и \ref{thBW1} из \S~\ref{S3} с некоторыми измениями в промежуточных утверждениях \ref{fgiii0}.  Условие \eqref{Lind},  очевидно, охватывает все  распределения вещественных точек ${\sf W}\subset \RR$ и распределения точек,  симметричные относительно вещественной  оси.
\end{enumerate}
Ключевую роль в доказательствах всех теорем из  \S~3 играет теорема \ref{th2L} из \S~\ref{SecL}  о взаимосвязях между 
логарифмическими (суб)мерами для распределений точек и различными вариантами условия Линделёфа \eqref{L} рода $1$ для распределений точек на $\CC$.
 
\begin{remark}\label{rem1}
Ясно, что условие  
${\sf W}\subset \CC_{\rh}$ можно в них заменить на расположение ${\sf W}\subset \CC_{\lh }$ целиком в левой полуплоскости с помощью зеркальной симметрии относительно мнимой оси. 

Анализ доказательств теорем A,  B и С показывает, что при доказательствах  импликаций   
\ref{fgii0}$\Longrightarrow$\ref{fgi0} каждой из теорем  нигде не используется расположение распределения точек ${\sf W}$ именно в правой полуплоскости $\CC_{\rh}$. Таким образом, условие    ${\sf W}\subset \CC_{\rh}$ в теоремах A,  B  и C применяется только при доказательстве импликаций \ref{fgi0}$\Longrightarrow$\ref{fgiii0}$\Longrightarrow$\ref{fgii0}. 
\end{remark}

\section{Условия Линделёфа рода $1$}\label{SecL}

Распределение точек 
${\sf Z}=\{{\sf z}_{\tt j}\}_{\tt j\in J}$  удовлетворяет {\it условию Линделёфа\/}  (рода $1$), если  
\begin{subequations}\label{L}
\begin{align}
\sup_{r\geq 1}\biggl|\sum_{1< |{\sf z}_{\tt j}|\leq r} \frac{1}{{\sf z}_{\tt j}}
\biggr|<+\infty; 
\tag{\ref{L}L}\label{con:LpZ} 
\\
\intertext{удовлетворяет {\it $\RR$-условию Линделёфа\/} (рода $1$), если}
\sup_{r\geq 1}\biggl| \sum_{1< |{\sf z}_{\tt j}|\leq r} \Re\frac{1}{{\sf z}_{\tt j}}
\biggr|<+\infty; 
\tag{\ref{L}R}\label{con:LpZR}\\
\intertext{удовлетворяет {\it $i\RR$-условию Линделёфа\/} (рода $1$), если (cм. \eqref{Lind})}
\sup_{r\geq 1}\biggl|\sum_{1< |{\sf z}_{\tt j}|\leq r} \Im \frac{1}{{\sf z}_{\tt j}}
\biggr|<+\infty. 
\tag{\ref{L}I}\label{con:LpZiR} 
\end{align}  
\end{subequations} 
Особую роль условия Линделёфа \eqref{con:LpZ}  рода $1$ обеспечивает   
\begin{thHL}[{\cite{Levin56}, \cite{Levin96}, \cite[2.10]{Boas}}]\label{pr:repef}
Если   $f\neq 0$ --- ц.ф.э.т., то распределение корней $\Zero_f$ конечной верхней плотности  
и удовлетворяет условию Линделёфа \eqref{con:LpZ}. Обратно, если  распределение точек ${\sf Z}$ конечной верхней плотности и  удовлетворяет условию Линделёфа \eqref{con:LpZ}, то найдётся  ц.ф.э.т. $f\neq 0$  с  $\Zero_f={\sf Z}$.
\end{thHL}

Следующее предложение 
сразу следуют  из определений  \eqref{L} и \eqref{logZC}--\eqref{df:l}.

\begin{proposition}\label{remZL}
Пусть ${\sf Z}\subset \CC$ --- распределение комплексных точек. 

\begin{enumerate}[{\rm [L1]}]
\item\label{L1}  ${\sf Z}$  удовлетворяет условию Линделёфа 
 \eqref{con:LpZ}, если  и только если   ${\sf Z}$ удовлетворяет одновременно  $\RR$-условию Линделёфа \eqref{con:LpZR} и 
$i\RR$-условию Линделёфа \eqref{con:LpZiR}.
\item\label{L2}  Следующие три утверждения эквивалентны:
\begin{enumerate}[{\rm (i)}]
\item\label{L2i} ${\sf Z}$ удовлетворяет  $\RR$-условию Линделёфа \eqref{con:LpZR}; 
\item\label{L2ii} для правых и левых логарифмических мер \eqref{df:dDl+} и \eqref{df:dDl-} выполнено соотношение
\begin{equation}\label{remLLrh}
\sup_{1\leq r<R<+\infty} 
\bigl|l_{\sf Z}^{\rh}(r,R)- l_{\sf Z}^{\lh}(r,R)\bigr|<+\infty;
\end{equation}
\item\label{L2iii} для логарифмической  субмеры  \eqref{df:dDlL} выполнено соотношение 
\begin{equation}\label{remLL}
\sup_{1\leq r<R<+\infty} 
\Bigl(\bigl|l_{\sf Z}^{\lh}(r,R)-l_{\sf Z}(r,R)\bigr|+ \bigl|l_{\sf Z}(r,R)- l_{\sf Z}^{\rh}(r,R)\bigr|\Bigr)<+\infty;
\end{equation}
\item\label{L2iv}  поворот $i{\sf Z}:=\{i{\sf z}_{\tt j}\}\subset \CC$ на угол $\frac{\pi}{2}$ удовлетворяет  $i\RR$-условию Линделёфа  \eqref{con:LpZiR}.
\end{enumerate}
\end{enumerate}
\end{proposition}


Для доказательства  последующих трёх теорем потребуется 
\begin{theorem}\label{th2L} Пусть ${\sf Z}\subset \CC$ --- распределение точек конечной верхней плотности. Тогда
существуют такие  
распределения вещественных точек ${\sf X}\subset \RR$  и ${\sf Y}\subset \RR$ 
конечной верхней плотности, что
объединение ${\sf Z}\cup {\sf X}$ удовлетворяет $\RR$-условию Линделёфа \eqref{con:LpZR}, 
объединение ${\sf Z}\cup {i\sf Y}$ удовлетворяет $i\RR$-условию Линделёфа \eqref{con:LpZiR},
 а также 
\begin{subequations}\label{XY}
\begin{align}
\sup_{1\leq r<R<+\infty}& 
\bigl(l_{{\sf Z}\cup {\sf X}}(r,R)- l_{\sf Z}(r,R)\bigr)<+\infty,
\tag{\ref{XY}X}\label{LZXY}
\\
\sup_{1\leq r<R<+\infty}& 
\bigl(l_{(i{\sf Z})\cup {\sf Y}}(r,R)- l_{i\sf Z}(r,R)\bigr)<+\infty.
\tag{\ref{XY}Y}\label{LZXYi}
\end{align}
\end{subequations}
При этом объединение ${\sf Z}\cup {\sf X}\cup i{\sf Y}$  удовлетворяет  условию Линделёфа  \eqref{con:LpZ}.
\end{theorem}  
\begin{proof} 
Рассмотрим двоичную последовательность $r_k:=2^k$ с $k\in \NN_0$. 
 
Для построения распределения точек ${\sf X}=\{{\sf x}_{\tt j}\}\subset \RR$ при каждом $k\in \NN_0$  произведём следующий выбор интервалов $I_k^{+}\subset \RR^+$ и $I_k^{-}\subset -\RR^+$, а также точек в них, образующих ${\sf X}$:
\begin{enumerate}[{\rm 1)}]
\item[{[+]}]\label{1X} Если $ l^{\rh}_{\sf Z}(r_k, r_{k+1})\leq l^{\lh}_{\sf Z}(r_k, r_{k+1})$, то положим  $I_k^+:=(r_k, r_{k+1}]\subset \RR^+$. Очевидно,  можно выбрать  $N_k\in \NN$ экземпляров точки $r_{k+1}\in I_k^+$ так, что  
\begin{multline}\label{chx+}
l_{\sf Z}(r_k, r_{k+1})\overset{\eqref{df:dDlL}}{=}l^{\lh}_{\sf Z}(r_k, r_{k+1})\leq 
l^{\rh}_{\sf Z}(r_k, r_{k+1})+N_k\frac{1}{r_{k+1}}
\\
\leq l^{\lh}_{\sf Z}(r_k, r_{k+1})+\frac{1}{r_{k+1}}\overset{\eqref{df:dDlL}}{=}  l_{\sf Z}(r_k, r_{k+1})+\frac{1}{r_{k+1}}.
\end{multline}
После этого отнесём эти $N_k$ точек $r_{k+1}$ на $I_k^+\subset \RR^+$ к распределению точек $X$.
\item[{[--]}]\label{2X} Если $l^{\lh}_{\sf Z}(r_k, r_{k+1})< l^{\rh}_{\sf Z}(r_k, r_{k+1})$, то  положим  $I_k^-:=[-r_{k+1} -r_k)\subset -\RR^+$.  Очевидно,  можно выбрать  $N_k\in \NN$ экземпляров точки $-r_{k+1}\in [-r_{k+1}, -r_k)=I_k^-$ так, что  
\begin{multline}\label{chx-}
l_{\sf Z}(r_k, r_{k+1})\overset{\eqref{df:dDlL}}{=}l^{\rh}_{\sf Z}(r_k, r_{k+1})\leq 
l^{\lh}_{\sf Z}(r_k, r_{k+1})+N_k\frac{1}{r_{k+1}}
\\
\leq l^{\rh}_{\sf Z}(r_k, r_{k+1})+\frac{1}{r_{k+1}} \overset{\eqref{df:dDlL}}{=} l_{\sf Z}(r_k, r_{k+1})+\frac{1}{r_{k+1}}.
\end{multline}
После этого также отнесём эти $N_k$ точек $-r_{k+1}$ на $I_k^-\subset -\RR^+$ к распределению точек $X$.
\end{enumerate}
Этим построение распределения точек ${\sf X}$ завершено. Построения  [$\pm$] показывают, что интервалы $I_k^+\subset \RR^+$ и $-I_k^-\subset \RR^+$ не пересекаются и покрывают весь луч $[1,+\infty)\subset \RR^+$.

Покажем сначала, что распределение точек ${\sf X}$ конечной верхней плотности. 

По построению из верхних оценок в \eqref{chx+} или в \eqref{chx-} с завершающим  равенством следует 
\begin{equation*}
N_k\frac{1}{r_{k+1}}
\leq l_{\sf Z}(r_k, r_{k+1})+\frac{1}{r_{k+1}}\overset{\eqref{df:divmn}}{\leq} 
\int_{\overline D(r_{k+1})\setminus \overline D(r_{k})} \Bigl|\Re \frac{1}{z}\Bigr| \dd n_{\sf Z}(z)+\frac{1}{r_{k+1}},
\end{equation*}
где $n_{\sf Z}$ --- считающая мера из \eqref{df:divmn}.  Переходя к радиальной считающей функции \eqref{nrad}, можем продолжить неравенства как 
\begin{equation}\label{Nk}
N_k\leq {r_{k+1}} \int_{r_k}^{r_{k+1}} \frac{1}{t}\dd n_{\sf Z}^{\rad}(t)+2\leq r_{k+1}\frac{1}{r_k}
\int_{r_k}^{r_{k+1}} \dd n_{\sf Z}^{\rad}(t) +1\leq 2\bigl(n_{\sf Z}^{\rad}(r_{k+1})- n_{\sf Z}^{\rad}(r_{k})\bigr)+1,
\end{equation}
откуда суммирование по $k$ даёт неравенства
\begin{equation*}
n_{\sf X}^{\rad}(r_{n+1})\leq \sum_{k=0}^{k=n} \Bigl(2\bigl(n_{\sf Z}^{\rad}(r_{k+1})- n_{\sf Z}^{\rad}(r_{k})\bigr)+ 1\Bigr) 
\leq 2n_{\sf Z}^{\rad}(r_{n+1})+(n+1) \quad\text{для всех $n\in \NN$}.  
\end{equation*}
Отсюда при каждом $n\in \NN$ для любого $r\in (r_n,r_{n+1}]$ получаем 
\begin{equation*}
\frac{n_{\sf X}^{\rad}(r)}{r}\leq 2\frac{n_{\sf X}^{\rad}(r_{n+1})}{r_{n+1}}
\leq 4\frac{n_{\sf Z}^{\rad}(r_{n+1})}{r_{n+1}}+2\frac{n+1}{2^{n+1}}=O(1) \quad \text{при $n\to \infty$ и $r\to +\infty$},
\end{equation*}
что доказывает конечность верхней плотности распределения точек ${\sf X}$. 

Покажем теперь, что объединение ${\sf Z}\cup {\sf X}$ удовлетворяет условию \eqref{LZXY}. 

При $n< N$ по построению  [+] части распределения точек ${\sf X}$, лежащей на положительной полуоси $\RR^+$, 
из аддитивности последовательно правой \eqref{df:dDl+} и левой  \eqref{df:dDl-} логарифмических мер и из  промежуточного срединного неравенства в \eqref{chx+}  следует
\begin{subequations}\label{lX}
\begin{align}
l_{{\sf Z}\cup {\sf X}}^{\rh}(r_n, r_N)&\overset{\eqref{df:dDl+}}{=} \sum_{k=n}^{N-1}l_{{\sf Z}\cup {\sf X}}^{\rh}(r_k, r_{k+1})\overset{\eqref{chx+}}{\leq}  \sum_{k=n}^{N-1}\Bigl(l_{{\sf Z}}^{\lh}(r_k, r_{k+1})+\frac{1}{r_{k+1}}\Bigr)
\notag
\\
&\leq \sum_{k=n}^{N-1}l_{{\sf Z}}^{\lh}(r_k, r_{k+1})+\sum_{k=n}^{N-1}\frac{1}{r_{k+1}}
\overset{\eqref{df:dDl-}}{=}l_{{\sf Z}}^{\lh}(r_n, r_N)+\sum_{k=n}^{N-1}\frac{1}{2^{k+1}}
\overset{\eqref{df:dDlL}}{\leq} l_{{\sf Z}}(r_n, r_N)+1.
\tag{\ref{lX}r}\label{lXZ-}\\
\intertext{Аналогично,  по построению  [--] части распределения точек ${\sf X}$, лежащей на отрицательной  полуоси $-\RR^+$, 
из аддитивности  последовательно  левой  \eqref{df:dDl-} и правой \eqref{df:dDl+}  логарифмических мер и из  промежуточного срединного  неравенства в \eqref{chx-}  следует}
l_{{\sf Z}\cup {\sf X}}^{\lh}(r_n, r_N)&\overset{\eqref{df:dDl-}}{=} \sum_{k=n}^{N-1}l_{{\sf Z}\cup {\sf X}}^{\lh}(r_k, r_{k+1})
\overset{\eqref{chx-}}{\leq}  \sum_{k=n}^{N-1}\Bigl(l_{{\sf Z}}^{\rh}(r_k, r_{k+1})+\frac{1}{r_{k+1}}\Bigr)
\notag
\\
&\leq \sum_{k=n}^{N-1}l_{{\sf Z}}^{\rh}(r_k, r_{k+1})+\sum_{k=n}^{N-1}\frac{1}{r_{k+1}}
\overset{\eqref{df:dDl+}}{=}l_{{\sf Z}}^{\rh}(r_n, r_N)+\sum_{k=n}^{N-1}\frac{1}{2^{k+1}}
\overset{\eqref{df:dDlL}}{\leq} l_{{\sf Z}}(r_n, r_N)+1.  
\tag{\ref{lX}l}\label{lXZ}
\end{align}
\end{subequations}

Из \eqref{lX} по определению  \eqref{df:dDlL} логарифмической субмеры имеем 
\begin{equation}\label{ll2}
l_{{\sf Z}\cup {\sf X}}(r_n, r_N)\overset{\eqref{df:dDlL}}{=}
\max \Bigl\{l_{{\sf Z}\cup {\sf X}}^{\lh}(r_n, r_N),  l_{{\sf Z}\cup {\sf X}}^{\rh}(r_n, r_N)\Bigr\}
\leq l_{{\sf Z}}(r_n, r_N)+1. 
\end{equation}
Отсюда  при любых $n\leq N$ для любых $r_n<r\leq r_{n+1}$ и $r_N<R\leq r_{N+1}$ при $r<R$ получаем 
\begin{multline*}
l_{{\sf Z}}(r, R)\leq l_{{\sf Z}\cup {\sf X}}(r, R)
\leq l_{{\sf Z}\cup {\sf X}}(r, r_{n+1})+
 l_{{\sf Z}\cup {\sf X}}(r_{n+1}, r_N)+ l_{{\sf Z}\cup {\sf X}}(r_N,r)\\
\overset{\eqref{ll2}}{\leq} 
\frac{1}{r_n}n_{{\sf Z}\cup {\sf X}}^{\rad}(r_{n+1})
+\bigl(l_{{\sf Z}}(r_{n+1}, r_N)+1\bigr)+
\frac{1}{r_N}n_{{\sf Z}\cup {\sf X}}^{\rad}(r_{N+1})\\
\leq 
l_{{\sf Z}}(r, R)+\Bigl(1 +2 \frac{n_{{\sf Z}\cup {\sf X}}^{\rad}(r_{n+1})}{r_{n+1}}
+2\frac{n_{{\sf Z}\cup {\sf X}}^{\rad}(r_{N+1})}{r_{N+1}}\Bigr).
\end{multline*}
Ввиду доказанной конечности верхней плотности распределения точек ${\sf X}$ последняя величина в круглой скобке ограничена по всем $0\leq n\leq N<+\infty$, что доказывает \eqref{LZXY}. 

Убедимся наконец, что объединение  ${\sf Z}\cup {\sf X}$ удовлетворяет $\RR$-условию Линделёфа \eqref{con:LpZR}. 
  
В случае [+] ввиду равенств 
\begin{equation*}
l^{\rh}_{\sf Z}(r_k, r_{k+1})+N_k\frac{1}{r_{k+1}}=l^{\rh}_{{\sf Z}\cup {\sf X}}(r_k, r_{k+1}) \quad\text{для всех  интервалов $(r_k, r_{k+1}]=I_k^+\subset \RR^+$}
\end{equation*}
из промежуточного срединного неравенства в    \eqref{chx+} получаем 
$$
0\leq l^{\rh}_{{\sf Z}\cup {\sf X}}(r_k, r_{k+1})-l_{\sf Z}^{\lh}(r_k, r_{k+1})\leq \dfrac{1}{r_{k+1}}=\dfrac{1}{2^{k+1}},
$$
а  $l_{\sf Z}^{\lh}(r_k, r_{k+1})=l_{{\sf Z}\cup {\sf X}}^{\lh}(r_k, r_{k+1})$, поскольку в противоположном интервале  $-I_k^+\subset -\RR^+$ по построению [+] точки для распределения точек   ${\sf X}$ не выбирались. Таким образом, 
\begin{equation}\label{lrhX}
0\leq l^{\rh}_{{\sf Z}\cup {\sf X}}(r_k, r_{k+1})-l_{{\sf Z}\cup {\sf X}}^{\lh}(r_k, r_{k+1})\leq \dfrac{1}{2^{k+1}}
\quad\text{для всех  интервалов $(r_k, r_{k+1}]=I_k^+\subset \RR^+$}.
\end{equation}

В случае [--] ввиду равенств 
\begin{equation*}
l^{\lh}_{\sf Z}(r_k, r_{k+1})+N_k\frac{1}{r_{k+1}}=l^{\lh}_{{\sf Z}\cup {\sf X}}(r_k, r_{k+1}), \quad\text{для всех  интервалов $(r_k, r_{k+1}]=-I_k^-\subset \RR^+$}
\end{equation*}
из промежуточного срединного неравенства в   \eqref{chx-} получаем 
$$
0\leq l^{\lh}_{{\sf Z}\cup {\sf X}}(r_k, r_{k+1})-l_{\sf Z}^{\rh}(r_k, r_{k+1})\leq \dfrac{1}{r_{k+1}}=\dfrac{1}{2^{k+1}},
$$
а  $l_{\sf Z}^{\rh}(r_k, r_{k+1})=l_{{\sf Z}\cup {\sf X}}^{\rh}(r_k, r_{k+1})$, поскольку в противоположном интервале  $-I_k^-\subset \RR^+$ по построению [--] точки для распределения точек  ${\sf X}$ не выбирались. Таким образом, 
\begin{equation}\label{lrhXl}
0\leq l^{\lh}_{{\sf Z}\cup {\sf X}}(r_k, r_{k+1})-l_{{\sf Z}\cup {\sf X}}^{\rh}(r_k, r_{k+1})\leq \dfrac{1}{r_{k+1}}\leq \dfrac{1}{2^{k+1}}
\quad\text{для всех  $(r_k, r_{k+1}]=-I_k^-\subset \RR^+$}.
\end{equation}
При  $r\in (r_n, r_{n+1}]$ суммирование по $k\leq n$ на основе 
\eqref{lrhX} и \eqref{lrhXl}
 даёт 
\begin{multline*}
\bigl|l^{\lh}_{{\sf Z}\cup {\sf X}}(1, r)-l_{{\sf Z}\cup {\sf X}}^{\rh}(1, r)\bigr|\leq
\sum_{k=0}^{n-1} \bigl|l^{\lh}_{{\sf Z}\cup {\sf X}}(r_k, r_{k+1})-l_{{\sf Z}\cup {\sf X}}^{\rh}(r_k, r_{k+1})\bigr|
+\bigl|l^{\lh}_{{\sf Z}\cup {\sf X}}(r_n, r)\bigr|\\
\overset{\eqref{lrhX},\eqref{lrhXl}}{\leq} \sum_{k=0}^{n-1}\frac1{2^{k+1}}+
\frac1{r_n}n_{{\sf Z}\cup {\sf X}}^{\rad}(r_{k+1})\leq 1+2
\frac{n_{{\sf Z}\cup {\sf X}}^{\rad}(r_{k+1})}{r_{n+1}}=O(1) \quad\text{при всех $n\in \NN$ и  $r\in (r_n, r_{n+1}]$.}
\end{multline*}
Последнее   из эквивалентности   (ii)$\Longleftrightarrow$(i)  части L\ref{L2} предложения \ref{remZL}
 с соотношением  \eqref{remLLrh} означает,  что для распределения точек ${\sf Z}\cup {\sf X}$ 
выполнено $\RR$-условие Линделёфа \eqref{con:LpZR}.

Существование  распределения точек ${\sf Y}$ с требуемыми свойствами, включая \eqref{LZXYi}, 
сразу  следует из  доказанной части теоремы \ref{th2L} для ${\sf X}$ после поворота $i{\sf Z}$ на угол $\frac{\pi}{2}$ распределения точек ${\sf Z}$ с учётом эквивалентности (i)$\Longleftrightarrow$(iv) части [L\ref{L2}] предложения \ref{remZL}. 
Условие Линделёфа \eqref{con:LpZ} для объединения ${\sf Z}\cup {\sf X}\cup i{\sf Y}$ следует из 
части [L\ref{L1}] предложения \ref{remZL}. 
\end{proof}
\begin{corollary}\label{corW}
Пусть распределение точек  ${\sf W}\subset \CC$ конечной верхней плотности. Тогда 
\begin{enumerate}[{\rm (a)}]
\item\label{Aa} Существует ц.ф.э.т. $g\neq 0$ со свойствами $g({\sf W})=0$ и 
\begin{equation}\label{supW}
\sup_{1\leq r<R<+\infty}\Bigl(l_{\Zero_g}(r,R)-l_{\sf W}(r,R)\Bigr)<+\infty.
\end{equation}
\item\label{Bb} Если\/ ${\sf W}$ асимптотически отделено углами от $i\RR$ в смысле \eqref{con:dis} 
и удовлетворяет $i\RR$-условию Линделёфа  \eqref{Lind}, то ц.ф.э.т. $g\neq 0$ со свойствами  ${\sf W}\subset \Zero_g$ и \eqref{supW} можно выбрать с асимптотически отделённым углами от $i\RR$ распределением корней  $\Zero_g$. 
\end{enumerate}
\end{corollary}
\begin{proof} По теореме \ref{th2L} существуют распределения точек ${\sf X}\subset \RR$
и $i{\sf Y}\subset i\RR$ конечной верхней плотности, для которых ${\sf W}\cup {\sf X}\cup i{\sf Y}$ удовлетворяет условию Линделёфа \eqref{con:LpZ} и выполнено \eqref{LZXY} с ${\sf W}$ вместо ${\sf Z}$. По теореме Адамара\,--\,Линделёфа существует ц.ф.э.т. $g\neq 0$ с распределением корней $\Zero_g={\sf W}\cup {\sf X}\cup i{\sf Y}$, обращающаяся в нуль на ${\sf W}$, для которой по \eqref{LZXY} с ${\sf W}$ вместо ${\sf Z}$ выполнено \eqref{supW}, поскольку распределение мнимых точек $i{\sf Y}\subset i\RR$ на  логарифмическую субмеру  никак не влияет. Часть \eqref{Aa} доказана.

По части \eqref{Bb} достаточно рассмотреть объединение ${\sf W}\cup {\sf X}$ с предыдущим выбором распределения вещественных точек ${\sf X}\subset \RR$. По $i\RR$-условию Линделёфа  \eqref{Lind} для ${\sf W}$ объединение  ${\sf W}\cup {\sf X}$  удовлетворяет условию Линделёфа \eqref{con:LpZ} по части [L\ref{L1}] предложения \ref{remZL}. По теореме Адамара\,--\,Линделёфа существует ц.ф.э.т. $g\neq 0$ с распределением корней $\Zero_g={\sf W}\cup {\sf X}$, обращающаяся в нуль на ${\sf W}$, для которой по \eqref{LZXY} с ${\sf W}$ вместо ${\sf Z}$ выполнено \eqref{supW}. Очевидно,  по-прежнему $\Zero_g={\sf W}\cup {\sf X}$ асимптотически отделено углами от $i\RR$, поскольку таковы как ${\sf W}$ по условию, так и распределение вещественных точек ${\sf X}\subset \RR$. 
\end{proof}

\section{Теоремы A, B и C без  условия расположения\\ распределения точек $\sf W$ в правой полуплоскости $\CC_{\rh}$}\label{S3}

Для целых функция $f\in \Hol_*(\CC)$ используем логарифмический интеграл 
\begin{equation}\label{fK:abp+}
J_{i\RR}(r,R;\ln |f|):=
\frac{1}{2\pi}\int_r^{R} \frac{\ln \bigl|f(iy)f(-iy)\bigr|}{y^2} \dd y, \quad 0<r<R<+\infty.
\end{equation} 

Вариант теоремы A с любым распределением точек  $\sf W$ конечной верхней плотности ---
\begin{theorem}\label{thAW}
Достаточно лишь  условия \eqref{dnsZW} конечности верхней плотности ${\sf Z}$ и ${\sf W}$, чтобы  каждое из утверждений\/ {\rm \ref{fgie}} и\/ {\rm \ref{fgiie}} из теоремы\/ {\rm A} было  эквивалентно утверждению 
\begin{enumerate}
\item[{\rm $\rm \widehat{II}$.}] 
При любом $\varepsilon \in \RR^+\setminus 0$ существует пара функций $f\in I_*^1({\sf Z})$ и  $g\in I^1_*({\sf W})$ cо свойством \eqref{supW},  а также  борелевское подмножество $E\subset \RR$, удовлетворяющие   \eqref{fgiRe}. 
\end{enumerate}
\end{theorem}
\begin{proof} По части \eqref{Aa} следствия  \ref{corW} существует  ц.ф.э.т. $g\neq 0$, с требуемыми в утверждении $\rm \widehat{II}$ свойствами, а по утверждению {\rm \ref{fgie}} теоремы A найдётся  $f\in I_*^1({\sf Z})$ , удовлетворяющая    \eqref{fgiRe},  что и доказывает импликацию  \ref{fgie}$\Longrightarrow$$\rm \widehat{II}$.

Если выполнено $\rm \widehat{II}$, то соответствующее интегрирование \eqref{fgiRe} даёт  
\begin{equation}\label{JJ}
J_{i\RR}\bigl(r,R;\ln |f|\bigr)\leq J_{i\RR}\bigl(r,R;\ln |g|\bigr)+\varepsilon \ln \frac{R}{r}
+C\quad \text{при всех\/ $1\leq  r<R<+\infty$}.
\end{equation}

\begin{lemma}[{\cite[(1.3)]{Kha89}, \cite[(0.4)]{kh91AA} и в явном  в \cite[предложение 4.1, (4.19)]{KhII}}]\label{lemJl} При любом фиксированном числе $r_0>0$ для любой ц.ф.э.т. $f\neq 0$ имеет место соотношение
\begin{equation}\label{Jll}
\sup_{r_0\leq r<R<+\infty}\max \Bigl\{\bigl|J_{i\RR}(r,R;\ln |f|)-l_{\Zero_f}^{\rh}(r,R)\bigr|,\;
\bigl|J_{i\RR}(r,R;\ln|f|)-l_{\Zero_f}^{\lh}(r,R)\bigr|\Bigr\}
<+\infty.
\end{equation}
\end{lemma}

 По лемме \ref{lemJl}, дважды применённой  соответственно к ц.ф.э.т. $f\neq 0$ и к ц.ф.э.т. $g\neq 0$, 
для некоторых чисел $C_f, C_g, M\in \RR^+$ получаем
\begin{multline*}
l_{\sf Z}(r,R)\overset{\eqref{df:dDlL}}{\leq} l_{\Zero_f}(r,R)\overset{\eqref{Jll}}{\leq} 
J_{i\RR}\bigl(r,R;\ln |f|\bigr) +C_f
\overset{\eqref{JJ}}{\leq} J_{i\RR}\bigl(r,R;\ln |g|\bigr)+\varepsilon \ln \frac{R}{r}+C+C_f
\\
\overset{\eqref{Jll}}{\leq}
l_{\Zero_g}(r,R)+C_g+\varepsilon \ln \frac{R}{r}
+C+C_f \overset{\eqref{supW}}{\leq}
l_{\sf W}(r,R)+M+C_g+\varepsilon \ln \frac{R}{r}
+C+C_f 
\end{multline*}
 при всех\/ $1\leq  r<R<+\infty$, где при последнем переходе использовано свойство-условие  \eqref{supW} из $\rm \widehat{II}$. Переобозначив здесь  сумму $M+C_g+C+C_f $ как число $C$, получаем соотношение  \eqref{Zlde}
утверждения {\rm \ref{fgiie}} из теоремы\/ {\rm A}. 

Справедливость импликации \ref{fgiie}$\Longrightarrow$\ref{fgie} из теоремы\/ {\rm A} уже отмечена в замечании \ref{rem1} 
для любых распределений точек ${\sf Z}$ и ${\sf W}$ конечной верхней плотности.
\end{proof}
Вариация теоремы B с $i\RR$-условием Линделёфа \eqref{Lind} на $\sf W$ вместо условия $W \subset \CC_{\rh}$ ---
\begin{theorem}\label{thBWB}
Пусть при \eqref{dnsZW}  и ${\sf Z}\subset \CC$, и  ${\sf W}$  асимптотически отделены углами от  $i\RR$ в смысле \eqref{con:dis},  а распределение точек ${\sf W}$ удовлетворяет $i\RR$-условием Линделёфа \eqref{Lind}. Тогда каждое из утверждений\/ {\rm \ref{fgi}} и\/ {\rm \ref{fgii}} из теоремы\/ {\rm B}  эквивалентно утверждению 
\begin{enumerate}
\item[{\rm $\rm \widehat{II}$.}] 
Существует пара $f\in I_*^1({\sf Z})$ и  $g\in I^1_*({\sf W})$ с  \eqref{supW}, удовлетворяющие   \eqref{fgiR}. 
\end{enumerate}
\end{theorem}
\begin{proof}
По части \eqref{Bb} следствия  \ref{corW} существует  ц.ф.э.т. $g\neq 0$, с требуемыми в утверждении $\rm \widehat{II}$ свойствами, а по утверждению {\rm \ref{fgi}} теоремы B найдётся  $f\in I_*^1({\sf Z})$ , удовлетворяющая    \eqref{fgiR},  что и доказывает импликацию  \ref{fgi}$\Longrightarrow$$\rm \widehat{II}$.

Если выполнено  утверждение $\rm \widehat{II}$, то интегрируя, как в \eqref{fK:abp+},  при каждом $r_0\in \RR^+\setminus 0$ для некоторой ограниченной функции $Q\colon [r_0,+\infty)\to \RR^+$,  удовлетворяющей  условию 
\begin{equation}\label{diy}
\lim_{x\to +\infty}\frac{Q(x)}{x}=0,
\end{equation}
и для некоторого числа $C_0\in \RR^+$ из  соотношения \eqref{fgiR}  получаем 
\begin{equation}\label{JJg}
J_{i\RR}\bigl(r,R;\ln |f|\bigr)\leq J_{i\RR}\bigl(r,R;\ln |g|\bigr)+
\int_r^R\frac{Q(y)}{y^2}\dd y+C_0\quad\text{при всех $0< r_0\leq r<R<+\infty$}.
\end{equation}   

\begin{lemma}[{\cite[следствие 2.1]{SalKha20}}]\label{corln}
Пусть $r_0\in \RR^+\setminus 0$. Если для  функции $Q\colon [r_0,+\infty) \to \RR^+$ выполнено \eqref{diy},
то найдётся  убывающая функция   $d\colon [r_0,+\infty)\to \RR^+$, для которой 
\begin{subequations}\label{dQ}
\begin{align}
\int_r^R\frac{Q(x)}{x^2}\dd t&\leq d(R)\ln\frac{R}{r}\quad \text{при всех $r_0\leq r<R<+\infty$},
\tag{\ref{dQ}Q}\label{{ad}A}\\
\lim_{R\to +\infty} d(R)&=0. 
\tag{\ref{dQ}d}\label{{ad}0}
\end{align}
\end{subequations} 
Если для функции $d\colon [r_0,+\infty) \to \RR^+$ выполнено \eqref{{ad}0}, то найдётся возрастающая функция $Q\colon [r_0,+\infty)\to \RR^+$, для которой выполнено \eqref{diy} и 
\begin{equation}\label{Qd}
d(R)\ln\frac{R}{r} \leq \int_r^R\frac{Q(x)}{x^2}\dd x
\quad \text{при всех $r_0\leq r<R<+\infty$}.
\end{equation}
\end{lemma}

По первой части леммы  \ref{corln}  неравенство \eqref{JJg} можно переписать как 
\begin{equation}\label{Jl}
J_{i\RR}\bigl(r,R;\ln |f|\bigr)\leq J_{i\RR}\bigl(r,R;\ln |g|\bigr)+
d(R)\ln \frac{R}{r}+C_0\quad\text{для всех $r_0\leq r<R<+\infty$},
\end{equation}   
где $d$ --- некоторая ограниченная функция со свойством \eqref{{ad}0}. 

 По лемме \ref{lemJl}, дважды применённой  соответственно к ц.ф.э.т. $f\neq 0$ и к ц.ф.э.т. $g\neq 0$, 
для некоторых чисел $C_f, C_g, M\in \RR^+$ получаем 
\begin{multline*}
l_{\sf Z}(r,R)\overset{\eqref{df:dDlL}}{\leq} l_{\Zero_f}(r,R)\overset{\eqref{Jll}}{\leq} 
J_{i\RR}\bigl(r,R;\ln |f|\bigr) +C_f
\overset{\eqref{Jl}}{\leq} J_{i\RR}\bigl(r,R;\ln |g|\bigr) +d(R)\ln \frac{R}{r}+C_0+C_f
\\
\overset{\eqref{Jll}}{\leq}
l_{\Zero_g}(r,R)+C_g +d(R)\ln \frac{R}{r}+C_0+C_f \overset{\eqref{supW}}{\leq}
l_{\sf W}(r,R)+M+C_g+d(R)\ln \frac{R}{r}+C_0+C_f 
\end{multline*}
 при всех\/ $r_0\leq  r<R<+\infty$, где использовано свойство-условие  \eqref{supW} из $\rm \widehat{II}$. Полагая  $C:=M+C_g+C_0+C_f $, получаем соотношение  \eqref{Zld}
утверждения {\rm \ref{fgii}} из теоремы\/ {\rm B}. 

Справедливость импликации \ref{fgii}$\Longrightarrow$\ref{fgi} из теоремы\/ {\rm B} уже отмечена в замечании \ref{rem1} 
для любых распределений точек ${\sf Z}$ и ${\sf W}$ конечной верхней плотности.
\end{proof}

Вариация теоремы C с $i\RR$-условием Линделёфа \eqref{Zld0} на $\sf W$ вместо условия ${\sf W}\subset \CC_{\rh}$  ---
\begin{theorem}\label{thBW1}
Пусть выполнены условия  теоремы\/ {\rm \ref{thBWB}.} Тогда каждое из утверждений\/ {\rm \ref{fgi0}} и\/ {\rm \ref{fgii0}}  из теоремы\/ {\rm C}  эквивалентно утверждению 
\begin{enumerate}
\item[{\rm $\rm \widehat{II}$.}] 
Существует пара $f\in I_*^1({\sf Z})$ и  $g\in I^1_*({\sf W})$ с  \eqref{supW}, удовлетворяющая    \eqref{fgiR0}. 
\end{enumerate}
\end{theorem}
\begin{proof} По части \eqref{Bb} следствия  \ref{corW} существует  ц.ф.э.т. $g\neq 0$, с требуемыми в утверждении $\rm \widehat{II}$ свойствами, а по утверждению {\rm \ref{fgi0}} из теоремы Мальявена\,--\,Рубела найдётся  $f\in I_*^1({\sf Z})$ , удовлетворяющая    \eqref{fgiR0},  что и доказывает импликацию  \ref{fgi0}$\Longrightarrow$$\rm \widehat{II}$.

Если выполнено $\rm \widehat{II}$, то логарифмирование и интегрирование неравенства \eqref{fgiR0} по отрезку $[r,R]$ с делением на $y^2$, как в  \eqref{fK:abp+}, даёт неравенства 
\begin{equation}\label{Jfg}
J_{i\RR}\bigl(r,R;\ln |f|\bigr)\leq J_{i\RR}\bigl(r,R;\ln |g|\bigr)\quad\text{для всех $r_0\leq r<R<+\infty$}.
\end{equation} 
По лемме \ref{lemJl}, применённой  к ц.ф.э.т. $f\neq 0$, 
для некоторого числа $C_f\in \RR^+$ получаем
\begin{equation*}
l_{\sf Z}(r,R)\overset{\eqref{df:dDlL}}{\leq} l_{\Zero_f}(r,R)\overset{\eqref{Jll}}{\leq} 
J_{i\RR}\bigl(r,R;\ln |f|\bigr) +C_f\overset{\eqref{Jfg}}{\leq}
 J_{i\RR}\bigl(r,R;\ln |g|\bigr)+C_f
\end{equation*}
для всех $r_0\leq r<R<+\infty$. Вновь по  лемме \ref{lemJl}, но применённой  уже к ц.ф.э.т. $g\neq 0$,  можем продолжить эту цепочку неравенств с некоторым числом $C_g\in \RR^+$ как  
\begin{equation*}
 l_{\sf Z}(r,R)\leq  J_{i\RR}\bigl(r,R;\ln |g|\bigr)+C_f\overset{\eqref{Jll}}{\leq} 
l_{\Zero_g}^{\rh}(r,R)+C_g+C_f\text{ при всех\/ $r_0\leq  r<R<+\infty$.}
\end{equation*}
Отсюда по условию \eqref{supW} из утверждения $\rm \widehat{II}$ для  некоторого числа $C_0\in \RR^+$ получаем  
$$
 l_{\sf Z}(r,R)\leq  l_{\Zero_g}(r,R)+C_g+C_f
\overset{\eqref{supW}}{\leq} l_{\sf W}(r,R)+C_0+C_g+C_f
\quad\text{при всех\/ $r_0\leq  r<R<+\infty$,}
$$
откуда при $C:=C_0+C_g+C_f\in \RR^+$ получаем утверждение  \ref{fgii0} теоремы C с  \eqref{Zld0}.

Справедливость импликации \ref{fgii0}$\Longrightarrow$\ref{fgi0}  уже отмечена в замечании \ref{rem1}.
\end{proof}


\end{document}